\documentclass[12pt,a4paper]{article}
\usepackage[utf8]{inputenc}
\usepackage[T1]{fontenc}
\usepackage[english]{babel}
\usepackage{cite}
\usepackage{lmodern}
\usepackage{amsmath, amsfonts}
\usepackage{multirow}
\usepackage{algorithm,algorithmicx,algpseudocode}
\usepackage{graphicx, placeins, caption, subcaption}

\graphicspath{{img/}}

\def\cvar{\texttt{CVaR}}
\def\var{\texttt{VaR}}
\def\E{\mathbb{E}}

\title{A Multicut Approach to Compute Upper Bounds for Risk-Averse SDDP}

\date{\today}
\author{Joaquim Dias Garcia\footnote{PSR, Brazil. joaquim@psr-inc.com} \and Iago Leal de Freitas\footnote{Brazil. iago.lealf@gmail.com} \and Raphael Chabar\footnote{PSR, Brazil. chabar@psr-inc.com} \and Mario V. F. Pereira\footnote{PSR, Brazil. mario@psr-inc.com}}

\begin{document}

\maketitle

\abstract{%

Stochastic Dual Dynamic Programming (SDDP) is a widely used and fundamental algorithm
for solving multistage stochastic optimization problems.
Although SDDP has been frequently applied to solve risk-averse models with the
Conditional Value-at-Risk (CVaR), it is known that the estimation of upper bounds
is a methodological challenge, and many methods are computationally intensive.
In practice, this leaves most SDDP implementations without a practical and clear stopping criterion.
In this paper, we propose using the information already contained in a multicut formulation of SDDP
to solve this problem with a simple and computationally efficient methodology.

The multicut version of SDDP, in contrast with the typical average cut,
preserves the information about which scenarios give rise to the worst costs,
thus contributing to the CVaR value.
We use this fact to modify the standard sampling method on the forward step so the average of
multiple paths approximates the nested CVaR cost. We highlight that minimal changes are
required in the SDDP algorithm and there is no additional computational burden for
a fixed number of iterations.

We present multiple case studies to empirically demonstrate the effectiveness of the method.
First, we use a small hydrothermal dispatch test case, in which we can write the
deterministic equivalent of the entire scenario tree to show that the method perfectly computes
the correct objective values.
Then, we present results using a standard approximation of the Brazilian operation problem
and a real hydrothermal dispatch case based on data from Colombia.
Our numerical experiments showed that this method
consistently calculates upper bounds higher than lower bounds for those risk-averse problems
and that lower bounds are improved thanks to the better exploration of the scenarios tree.
}

\section{Introduction}

There are many problems that can be modeled using the framework of
multistage stochastic programming, especially in the energy sector.
For treating large problems,
the main algorithm in use today is called Stochastic Dual Dynamic Programming (SDDP),
introduced in~\cite{sddp}.
It is a clever variation of traditional stochastic dynamic programming
that samples paths (or trajectories) in a extremely large scenario tree and at each iteration uses linear programming duality
to approximate the future cost functions, also known as cost-to-go functions, for each time step.
When the problems considered are all linear,
the algorithm returns both a lower bound and an upper bound
that can be used to guarantee convergence to the optimal solution.

Traditionally,
stochastic dynamic programming represents the relations between stages
through a recursive equation known as the Bellman equation,
where we optimize the expected value among all possible realizations
of the random process for the next stage:
\begin{equation}
\label{eq:bellman_e}
   \begin{array}[t]{rl}
   Q_t(x_{t-1}, \xi_t) = \\
    \min\limits_{x_t,u_t} & C_t(x_{t-1}, \xi_t, u_t)+ \E_{\xi_{t+1}}[Q_{t+1}(x_{t}, \xi_{t+1})] \\
    \textrm{s.t.} & x_t = D_t(x_{t-1}, \xi_t, u_t) \\
                  & u_t \in U_t(x_{t-1}, \xi_t) 
  \end{array} 
\end{equation}
where
$x_{t-1}$ is the initial state and $\xi_t$ is the current value of the random variable,
both considered input coefficients for the optimization problem.
The decision variables are the state at the end of the stage $x_t$ and the control (or action) $u_t$.
$U_t$ is the set of feasible actions, $D_t$ is the state transition dynamics,
and $C_t$ is the immediate cost.
For the final stage $T$, we set $Q_T(\cdot, \cdot) = 0$.
This is called the \emph{risk-neutral} formulation of a multistage stochastic program.

In many practical problems, it is desirable to protect oneself
against the worst outcomes of random processes.
In this setting, called \emph{risk-averse}, the expected value, $\E_{\xi_{t+1}}[\cdot]$, is replaced
by a risk measure, $\rho_{\xi_{t+1}}[\cdot]$, typically a coherent risk measure~\cite{coherentMR}.
The most common of these risk measures is the \emph{Conditional Value-at-Risk} (CVaR)~\cite{cvar_rockafellar2000,cvar_fabian2008},
which is roughly equivalent to taking the average conditioned to be above a given quantile. The CVaR is also known as Average Value-at-Risk, AVaR, and Expected Shortfall, ES.
Unfortunately, in the presence of risk aversion,
we lose the capacity to properly estimate an upper bound for the SDDP~\cite{philpott2011}.

There are many recent advancements in the literature
that attempt to address the upper-bound problem through varying methods,
such as statistical samplings based on bad outcomes~\cite{guigues2021riskaverse,shapiro2020},
inner approximations of the Bellman function~\cite{philpott2013},
or dual variants of SDDP~\cite{leclere2020,dacosta2021dual}.

In this work, we present a method for estimating the upper bound
of risk-averse SDDP by modifying the sampling of scenarios in the forward pass through the weights endogenously defined by the optimal solution of the CVaR formulation of the future cost function. 
This method necessarily uses a multicut approximation for the future cost function and only requires the solution of the subproblem of the current stage to sample the random scenario of the next stage in the forward pass.
The method does not require any additional modification on SDDP's backward step
nor keeping annotations of sampled outcomes throughout iterations, which makes the method especially simple to implement.

In Section 2, we will detail the methodology briefly described above. Section 3 contains case studies with three test systems of hydrothermal dispatch: first, a small case that can be solved by a deterministic equivalent, then a standard academic representation of the Brazilian system, and finally, a realistic representation of the Colombian system. Section 4 outlines conclusions and future work.

\section{Methodology}


Throughout this paper,
we consider the following definition of $\cvar_\alpha[Y]$ of a random variable $Y$ with distribution $F_Y(y) = P(Y\leq y)$:
\[
    \cvar_\alpha[Y] = \E[Y | Y \geq \var_{\alpha}[y]] = \min_b \left\{b + \frac{1}{1-\alpha} \E[(Y - b)^+]\right\}
\]
where
\[
    \var_\alpha[Y] = \min_r \left\{F_Y(r) \geq \alpha  \right\}
\]
Also, we consider risk-averse SDDP
using a convex combination between the expected value and CVaR as a risk measure:
\[
  \rho[Y] = (1-\lambda)\E[Y] + \lambda\cvar_\alpha[Y]
\]
where $\lambda, \alpha \in [0, 1]$ are fixed parameters. Special cases are $\lambda = 0$
in which only the expected value is considered, $\lambda = ,1$ in which only the CVaR is considered,
$\alpha = 0$ in which the CVaR degenerates to the expected value,
and $\alpha = 1$ in which the CVaR represents the worst-case scenario as in a robust optimization framework.

This is a coherent risk measure~\cite{coherentMR} and can be
represented as the optimal solution to a linear program~\cite{cvar_rockafellar2000,cvar_fabian2008}.
For the case of a discrete random variable $Y$ with possible values $\{y^l\}_{l=1}^L$
all with probability $\frac{1}{L}$, we can write:
\begin{equation}
  \label{eq:lp_cvar}
  \rho[Y] = \begin{array}[t]{rl}
    \min\limits_{\delta^l, z} & \frac{1-\lambda}{L}\sum_{l=1}^L y^l + \lambda z  + \frac{\lambda}{(1-\alpha) L} \sum_{l=1}^L \delta^l \\
    \textrm{s.t.} & \delta^{l} \ge y^l - z, \quad \forall l =1, ...,L \\
                  & \delta^l \ge 0, \quad \forall l = 1, ..., L
  \end{array}
\end{equation}

We consider a multi-stage stochastic program whose cost function
satisfies the recursion:
\begin{equation}
\label{eq:bellman}
   \begin{array}[t]{rl}
   Q_t(x_{t-1}, \xi_t) = \\
    \min\limits_{x_t,u_t} & C_t(x_{t-1}, \xi_t, u_t)+ \rho_{\xi_{t+1}}[Q_{t+1}(x_{t}, \xi_{t+1})] \\
    \textrm{s.t.} & x_t = D_t(x_{t-1}, \xi_t, u_t) \\
                  & u_t \in U_t(x_{t-1}, \xi_t) 
  \end{array} 
\end{equation}
Which is the risk-averse version of the model \eqref{eq:bellman_e}.



Combining \eqref{eq:lp_cvar} and \eqref{eq:bellman},
and assuming that the possible outcomes of $\xi_{t+1}$ are $\{\xi_{t+1}^l\}_{l=1}^L$,
all with probability $\frac{1}{L}$, we get the nested CVaR multicut reformulation:
\begin{equation}
  \label{eq:recursion_formula}
  \begin{array}[t]{rl}
  Q_t(x_{t-1}, \xi_t) = \\
    \min\limits_{x_t,u_t,\delta_t^l, z_t, \beta_t^l } & C_t(x_{t-1}, \xi_t, u_t) + \frac{1-\lambda}{L} \sum_{l=1}^L \beta_t^l +\\
                   &\hspace{29mm} \lambda z_{t}  + \frac{\lambda}{(1-\alpha) L} \sum_{l=1}^L \delta_{t}^l \\
    \textrm{s.t.} & x_t = D_t(x_{t-1}, \xi_t, u_t) \\
                  & u_t \in U_t(x_{t-1}, \xi_t) \\
                  & \beta_t^l \ge Q_{t+1}(x_{t}, \xi_{t+1}^l), \quad \forall l = 1, ..., L \\
                  & \delta_{t}^l \ge \beta_t^l - z_{t}, \quad \forall l = 1, ..., L \\
                  & \delta_t^l \ge 0, \quad \forall l = 1, ..., L
  \end{array}
\end{equation}
Note that there is one explicit future cost function for each possible scenario of the next stage $\{\xi_{t+1}^l\}_{l=1}^L$.



\subsection{Stochastic dual dynamic programming}

Multi-stage stochastic programming
consists of a sequence of decision processes,
each one depending on the preceding and on the realization of random variables.
In theory, the random process can be continuous so that infinitely many outcomes are possible in each stage. However, for the application of the SDDP algorithm, it is common to consider a discretized version of the problem, in which all possible realizations of a random process can be presented in the form of a scenario tree.
This discretization step is the first application of Monte Carlo sampling required by SDDP practitioners.

\begin{figure}[h!]
\begin{center}
\subfloat[Original tree]{
\includegraphics[height=1.3in]{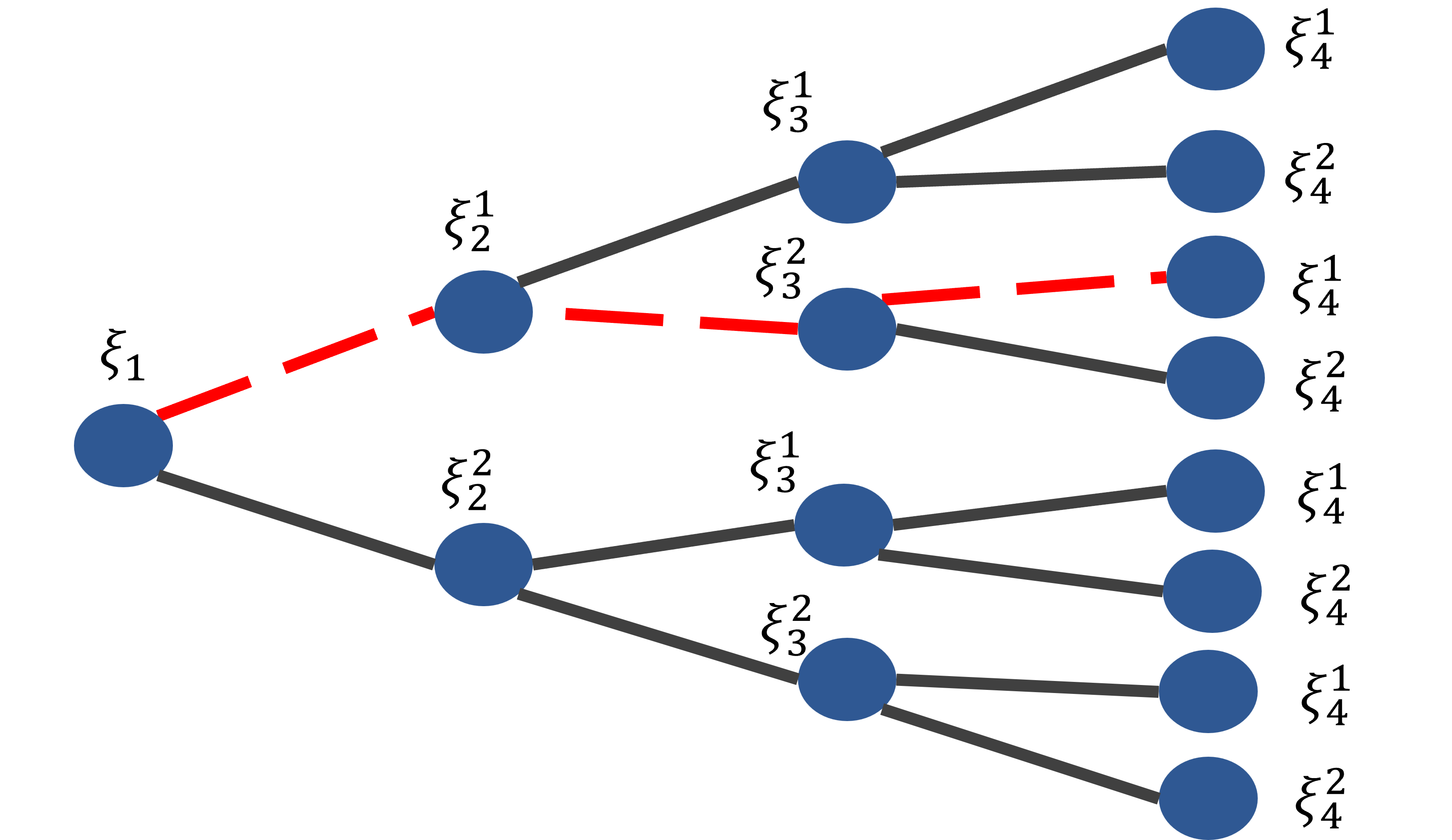}
}
\subfloat[Recombining tree]{
\includegraphics[height=1.3in]{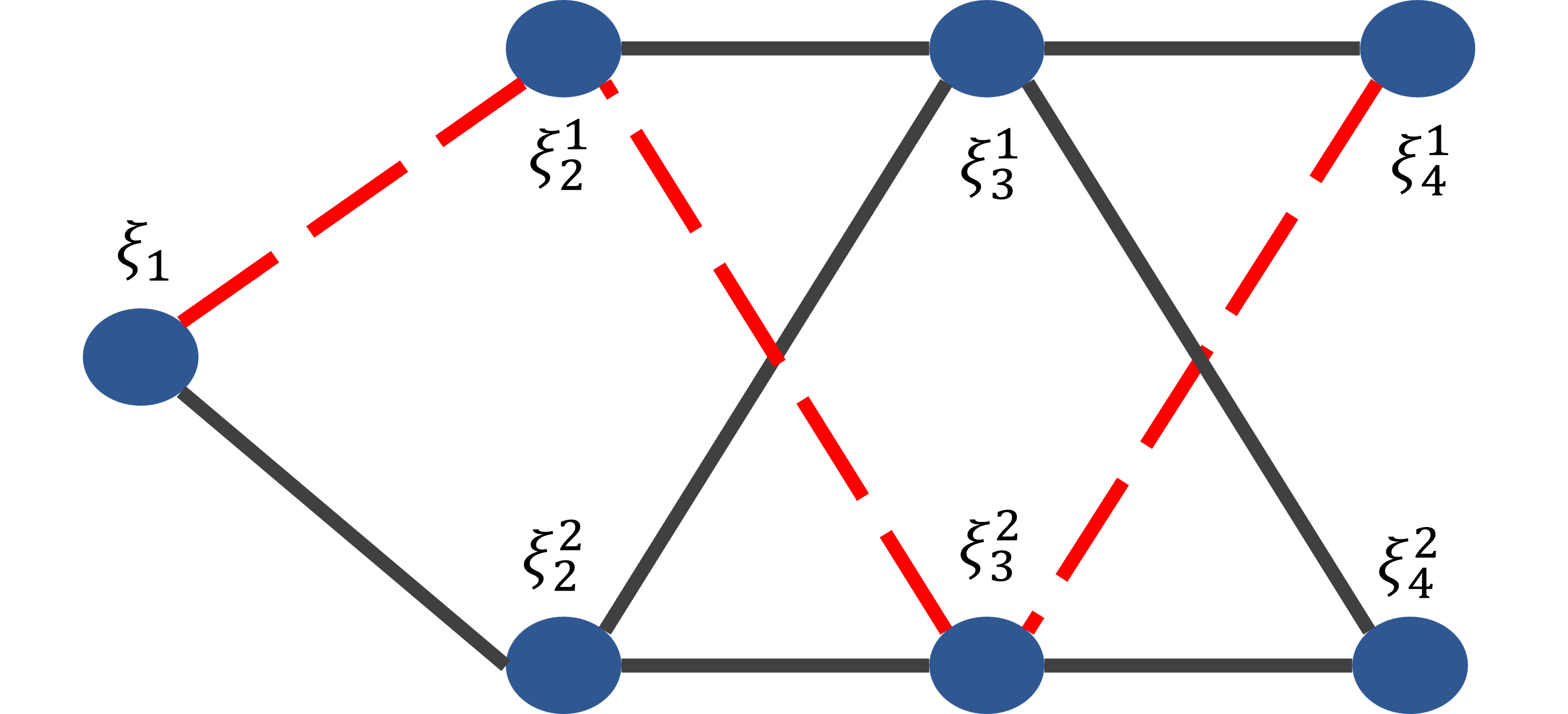}
}
\end{center}
\caption{Equivalence between a stagewise independent scenario tree and a recombining scenario tree. The dashed red trajectory is equivalent in both representations.} \label{fig:tree}
\end{figure}

To apply the SDDP algorithm, we make an additional requirement: random variables must be \textit{stagewise independent}. That is, the uncertainty $\xi_t$ is independent of all the previous stages $1$ to $t-1$.
The root node is the deterministic first-stage variable $\xi_1$
and for each node at stage $t$, there are the same possible realizations of the random process in the next stage, $t+1$:  $\{\xi_{t+1}^l\}_{l=1}^L$
In this particular case, the scenario tree degenerates into the so-called \textit{recombining scenario tree}, which is not an actual tree, but a graph that compactly represents a scenario tree with stagewise independency. Figure \ref{fig:tree} presents a stagewise independent scenario tree and the corresponding recombining scenario tree, a path from the first to the last stage is highlighted in red as an example.

The $Q$ functions are not known, in fact, they are what is being approximated. In SDDP, these future cost functions are approximated as the maximum of cutting planes. We name these approximations as $\tilde{Q}(x)$. In more detail, we define:
\begin{equation}
\tilde{Q}(x) = \max \{\pi_i^{\top}(x-\bar{x}_i)+q_i | i \in \mathcal{C}\},\label{eq:fcf}
\end{equation}
where the triplets $(\pi_i,\bar{x}_i, q_i)$ are the cut coefficients and $\mathcal{C}$ is the set of indices of cuts. At each iteration of the SDDP algorithm, more cuts will be added, improving the approximation. By construction $\tilde{Q}(x) \le {Q}(x), \forall x$, hence an approximation from below. In the first appearances of SDDP in the literature $\rho[Q]$ was approximated as a single function in the so-called single-cut version. Here, we apply the multicut representation, in which $\rho[Q]$ is first explicitly expanded as a combination of $L$ functions, one of each of the possible uncertainty realizations in the next stage, as described above, and then, each of these functions is approximated individually by cutting planes. Applying the  multicut approximation to \eqref{eq:recursion_formula} leads to:

\begin{equation}
  \label{eq:multicut}
  \begin{array}[t]{rl}
  \tilde{Q}_t(x_{t-1}, \xi_t) = \\
    \min\limits_{x_t,u_t,\delta_t^l, z_t, \beta_t^l } & C_t(x_{t-1}, \xi_t, u_t) + \frac{1-\lambda}{L} \sum_{l=1}^L \beta_t^l +\\
                   &\hspace{29mm} \lambda z_{t}  + \frac{\lambda}{(1-\alpha) L} \sum_{l=1}^L \delta_{t}^l \\
    \textrm{s.t.} \quad x_t &= D_t(z_t, \xi_t, u_t) \\
                  u_t &\in U_t(z_{t}, \xi_t) \\
                  z_t &= x_{t-1} \quad \quad : \ \pi_{t-1} \\
                  \beta_t^l &\ge \pi_{t,l,i}^\top (x_t-\bar{x}_{t,l,i}) - q_{t,l,i}, \forall  i \in \mathcal{C}_t^l, l = 1, ..., L \\
                  \delta_{t}^l &\ge \beta_t^l - z_{t}, \quad \forall l = 1, ..., L \\
                  \delta_t^l &\ge 0, \quad \forall l = 1, ..., L
  \end{array}
\end{equation}
Where we added the equation $z_t = x_{t-1}$ that simplifies the construction of the cut,
which will only need the dual of the constraint, $\pi_{t-1}$, the objective value $q_{t-1}$ and the input state $x_{t-1}$.
We denote the immediate cost, the evaluation of $C_t$ in the optimal solution, as $c_t$.

This formulation leads to larger optimization problems, but more detailed information is passed along. Also, the multicut approach is expected to be especially effective if the number of scenarios is significantly smaller than the number of constraints in the subproblem \cite{birge2011introduction}. Moreover, the explicit representation of the FCF of each node in the next stage will be instrumental to the next developments.

We can summarize a standard version of the SDDP algorithm as follows. Advanced versions of the algorithm, including Markov representations of uncertainties, infinite horizons, hazard-decision formulations, non-homogeneous trees, and so on, are possible, see \cite{dowson2020policy, street2020assessing}.

We start with the \textit{forward pass}. Following its name, the forward pass starts from the first node, with initial state $x_0$ and (deterministic) uncertainty realization $\xi_1$, and problems are solved in the forward direction of time. For each stage, $t$, only one node problem is solved given its past, the current realization of the random variable, $\xi_{t}^l$, and an uncertain view of the future. This uncertain view must ensure causality and non-anticipativity, thus, the future is not known and the decision must be taken assuming that all nodes in the next stage are equally probable. Hence, problem \eqref{eq:recursion_formula} is solved considering the future cost functions for each $l \in 1,..., L$.
After such a problem is solved, we must sample the next realization of the uncertainty $\xi_{t+1}^l$ so that the next problem can be solved until the process ends at stage $T$. A uniform sample is performed in the traditional SDDP. The sampling step is the second application of Monte Carlo and it will be the key to our proposed method.

The average cost of multiple of these solution paths is used to estimate upper bounds. The solution to the problem in the first stage leads to a deterministic lower bound, as opposed to the estimated statistical upper bound. Many stopping rules have been proposed after the first one by \cite{sddp}. In this work, any rule based on the deterministic lower bound and the estimated statistical upper bound can be considered.

After one or multiple scenario path is solved, we move to the \textit{backward pass}. A scenario path, $s$, solution is fully characterized by the set of tuples $\{(x_{t-1}^s, \xi_{t}^s)\}_{t\in \{1, ... , T\}}$ of optimized states and sampled uncertainties. As the name hints, the backward pass will solve problems in the reverse direction of time. Starting at stage $T$, all nodes associated with this stage will be solved, but all are conditioned to the input coming from the solutions path: $(x_{T-1}^s, \xi_{T}^s)$. The solutions of each node $l \in L$ will lead to a cut, $i$, containing $(\pi_{t,l,i}, q_{t,l,i})$, that will be stored in a cut pool indexed by $\mathcal{C}_t^l$ and associated with that node. The stagewise independency combined with the cut construction (a function of state and uncertainty) allows cut sharing so that the cuts generated by the current path are valid for all other possible paths. The process continues up to the second stage, which generates cuts for the first stage.

We summarize the SDDP algorithm in Algorithm~\ref{alg:sddp}. We present a batched version of the algorithm. Instead of creating a single solutions path, $s$, per iteration, we create $S$ scenario paths. This is important because we measure iterations as batches of scenarios. Also, this highlights the actual implementation used that solves each scenario path in parallel.
We also fix a maximum number of iterations, $K$, so that the algorithm stops if convergence is too slow.
For further details, see~\cite{sddp}.

\begin{algorithm}
    \caption{Multicut Batched Traditional SDDP}
    \label{alg:sddp}
    \begin{algorithmic}[1]
      \For {$k = 1, ..., K$}
        \State{(Forward Step)}
          \For{$t = 1, ..., T$}
        \For{$s = 1, ..., S$}
            \State Sample $\xi_{t}^s$ from $\{\xi_{t}^l\}_{l=1}^L$
            \State Solve \eqref{eq:recursion_formula} to obtain the state $x_{t}^s$ and immediate cost $c_t^s$ 
          \EndFor
        \EndFor
        \State{(Convergence Check)}
        \State Estimate statistical upper bound from $\{c_t^s\}_{s = 1, ..., S, t = 1, ..., T }$
        \State Store deterministic lower bound $c_1$
        \If {Stopping condition is met}
            \State Stop
        \EndIf
        \State{(Backward Step)}
        \For{$t = T, ..., 2$}
        \For{$s = 1, ..., S$}
          \For{$l = 1, ..., L$}
            \State Solve \eqref{eq:recursion_formula} to obtain the cut $i$ represented by $\pi_{t,l,i}, \bar{x}_{t,l,i}, q_{t,l,i}$
            \State Improve the representation of \eqref{eq:fcf}
          \EndFor
        \EndFor
        \EndFor
      \EndFor
    \end{algorithmic}
\end{algorithm}

\subsection{Risk adjusted scenario sampling}

The unique feature of the multicut formulations is
that it uses a different variable $\beta^l$ for each scenario of the next stage.
Therefore the linear programming formulation of the CVaR provides us with enough information
to know which scenarios are actually considered
when computing the \cvar.

We know from the literature~\cite{cvar_rockafellar2000}
that at the optimal of equation~\eqref{eq:multicut},
the decision variable $z_t$ will be equal to the Value-at-Risk
of the variables $\beta_t^l$ for $l \in 1,...,L$.
Furthermore, with the exception of degenerate cases,
the decision variables $\delta_{t}^l$
are non-zero only if the scenario $l$ is used to compute the \cvar\ value,
they are equal to how much $\beta_t^l$ is above the \var.
Since the $\beta_t^l$ represent the cost-to-go for each opening,
we may view their weight on the objective as representing
the contribution of each scenario for calculating the cost.

To make the weights explicit,
we start ordering the scenarios by the value of $\beta_t^l$ such that
\[\beta_t^{(1)} \le ... \le \beta_t^{(l)} \le ... \le \beta_t^{(L)},\]
and call $\nu$ the last opening such that $\delta_t^{(l)} = 0$
(that is, the cost-to-go is above the Value-at-Risk).
From the previous discussion, we know that $\beta_t^\nu = z_t$.
Then, by ordering the scenario we get:
\[ \begin{aligned}
  \rho_{\xi_{t+1}}[Q_{t+1}(\cdot)] &= \frac{1-\lambda}{L} \sum_{l=1}^L \beta_t^l + \lambda z_{t}  + \frac{\lambda}{(1-\alpha) L} \sum_{l=1}^L \delta_{t}^l \\
          &= \frac{1-\lambda}{L} \sum_{l=1}^L \beta_t^{(l)} + \lambda z_{t}  + \frac{\lambda}{(1-\alpha) L} \sum_{l=\nu}^{L} (\beta_{t}^{(l)} - z_t) \\
          &= \frac{1-\lambda}{L} \sum_{l=1}^L \beta_t^{(l)} + \lambda \beta_t^\nu  + \frac{\lambda}{(1-\alpha) L} \sum_{l=\nu+1}^{L} (\beta_{t}^{(l)} - \beta_t^\nu) \\
          &= \sum_{l=1}^{\nu-1} \frac{1-\lambda}{L} \beta_t^{(l)} + \left(\frac{1-\lambda}{L} + \lambda - \frac{\lambda(L-\nu)}{(1-\alpha) L} \right) \beta_{t}^\nu \\
          & \hspace{27mm}+ \sum_{l=\nu+1}^{L} \left( \frac{1-\lambda}{L} + \frac{\lambda}{(1-\alpha) L} \right) \beta_{t}^{(l)}
\end{aligned} \]

From the above, we note that computing the risk measure
is the same as taking a linear combination of the cost-to-go at all openings
with weights
\[
  w_{t}^l = \begin{cases}
    \frac{1-\lambda}{L} + \frac{\lambda}{(1-\alpha) L} ,& \beta_t^l > z_t   \\
    \frac{1-\lambda}{L} + \lambda - \frac{\lambda(L-\nu)}{(1-\alpha) L},& \beta_t^l = z_t   \\
    \frac{1-\lambda}{L},& \beta_t^l < z_t. \\
  \end{cases}
\]
Notice also that, for any stage $t$, $\sum_{l=1}^L w_t^l = 1$.
Thus, the weights represent a probability distribution on the openings
such that calculating the risk measure is the same
as calculating the average with respect to the $w_t^l$.

Given the above-defined weights, it is only necessary to modify line 5 of Algorithm 1 to obtain a version of SDDP that is able to compute the upper bounds of the CVaR-based risk-averse multistage stochastic program.
Clearly, there is no performance degradation on the computational cost of each iteration and the changes required are minimal: obtain the weights and change the sampling method.

Our implementation extends SDDP by using
the optimal decision variables $\beta_t^l, \delta_t^l, z_t$
calculated at the backward step
to discover the weights $w_t^l$ of choosing opening $l$ at stage $t$.
This is then applied on the forward step
to sample the scenarios according to this new probability.

This method can be interpreted as an endogenous importance sampling scheme to
compute the weight distributions used in the forward pass to obtain scenario paths.
As the SDDP progressively improves the approximation of the future cost function,
the weights obtained at each iteration improve until they arrive at a distribution
that leads to the upper bound of the problem.

In the special case of pure CVaR as a risk measure ($\lambda = 1$), the weights on non-CVaR contributing scenarios go to zero, which might lead to issues in exploring the scenario tree. This issue might affect the convergence. In this case and in similar cases where $\lambda$ is close to zero, one possible alternative is to proceed as follows: in even iteration ($k$ is even), use the risk-adjusted sampling approach proposed here, in odd iteration ($k$ is odd), use standard sampling, but do not compute upper bound nor check convergence. Of course, other alternating schemes might be used to improve the exploration of the scenario tree.

\begin{figure}[t]
  \centering
  \includegraphics[width=0.6\linewidth]{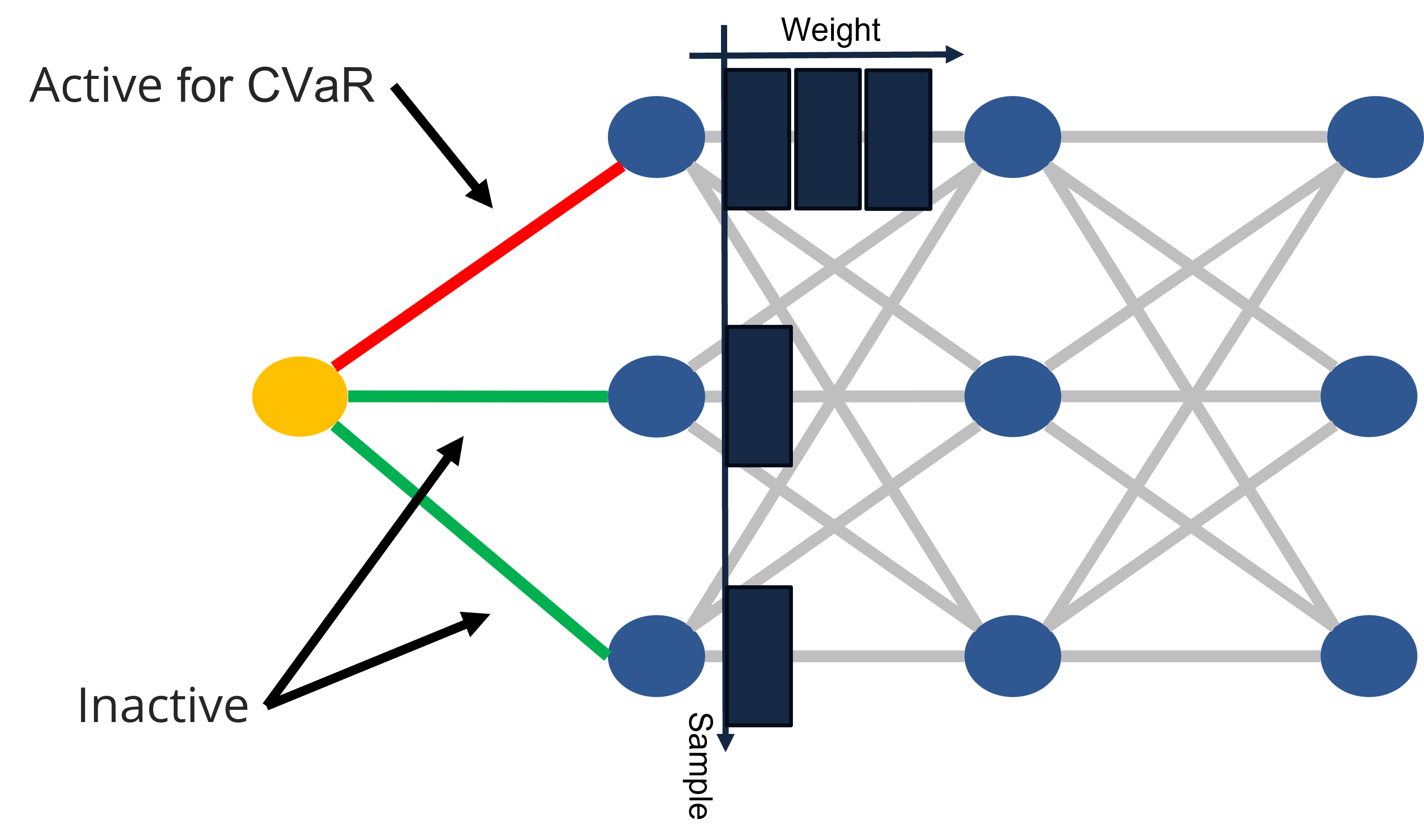}
  \caption{Illustration of sampling method for upper bound estimation of risk-averse SDDP.}
  \label{fig:cvar-sample}
\end{figure}

We finish with an illustration of the proposed method (with $0 <\lambda < 1$). Figure \ref{fig:cvar-sample} presents a depiction of the forward pass right after solving the subproblem of the first stage (marked in yellow). At this point, we know which scenarios contribute to the CVaR (line in red) and which do not (lines in green). Consequently, we can compute the weight that will be used to sample the scenario in stage 2, which is depicted in the plot of weights as a function of the sample that has more mass in the sample associated with the scenario active for CVaR. 

\section{Numerical experiments}

In this section, we compare the variants of SDDP using traditional sampling,
which leads to a naive upper bound known to be incorrect,
and with our modified sampling method based on multicut formulation,
which will lead to better upper bounds.
We apply the methods to the hydrothermal dispatch model that will be detailed next.
We consider 3 case studies:
first, a simple problem in which we can evaluate the entire tree, second,
a classical model for the Brazilian interconnected power system
found in~\cite{shapiro_case_br},
and, third, a case with real data from the Colombian system. Data for these cases are summarized in Table~\ref{tab:parameters}.
In all cases we used the \cvar\ parameters $\alpha = 0.5$
and $\lambda = 0.5$.

\begin{table}[b]
  \caption{Study parameters for numeric experiments.}
  \label{tab:parameters}
  \begin{center}
    \begin{tabular}{|c|c|c|c|}
        \hline
        Parameter & Full Tree & Brazil & Colombia \\
        \hline
        Hydro Reservoirs       & 2 & 4  & 38   \\
        Hydro Plants     & 2 & 4  & 136  \\
        Thermal Plants   & 3 & 95 & 202  \\
        Renewable Plants & 0 & 0  & 12   \\
        \hline
        Stages           & 10 or 7 & 120 & 60   \\
        Max autoregressive size      & 0 & 0   & 6    \\
        \hline
    \end{tabular}
  \end{center}
\end{table}

\subsection{Hydrothermal dispatch model}

In the following model, we have the sets:
$J^L$ is the set of lags in an autoregressive model,
$J^G$ is the set of thermal plants,
$J^R$ is the set of renewable plants,
$J^H$ is the set of hydro plants,
$J^F_n$ is the set of other nodes or areas connected to node $n$,
$J^U_j$ is the set of hydro plants upstream of plant $j$.
An additional $n$ as a subscript in $J^G$, $J^R$, and $J^H$ restricts the generators to the ones the node $n$.
The coefficients of the model are defined as:
$\tilde{D}_n(\xi)$ is the demand of subsystem (node) $n$ and scenario $\xi$,
$F_{n,m}$ is the maximum energy flow between nodes $n$ and $m$,
$C_j$ is the operating cost of thermal $j$,
$G_j$ is the maximum generation of thermal $j$,
$\tilde{R}_j(\xi)$ is the maximum generation of renewable $k$ at scenario $\xi$
$V_j$ is the maximum storage of hydro $j$,
$U_j$ is the maximum flow through the turbine of hydro $j$,
$p_j$ is the production coefficient of hydro $j$,
$\phi_{j,k}$ is the inflow auto-regressive coefficient of hydro $j$ and lag $k$,
$\tilde{\varepsilon}_i^t(\xi)$ is the inflow noise coefficient of hydro $j$ and scenario $\xi$.
Finally, the decision variables of the optimization problem are:
$g_j$, the generation of thermal $j$;
$r_j$, the generation of renewable $j$;
$f_{n,m}$, the flow of energy leaving node $n$ and reaching node $m$;
$u_j$, the turbine flow at hydro $j$;
$z_j$, the spill flow at hydro $j$;
$v_j^{t}$, the storage at hydro $j$, at the beginning of stage $t$, and at the end of stage $t-1$;
$a_j^t$, the inflow at hydro $j$, stage $t$.
$a^{[t]}$ stands for the vector of all inflows before stage $t$.

\begin{align}
  & {Q}_{t}\big(\{{v}^{t}_j, {a}^{[t-1]}_j\}_{j\in J^H}, \xi_{t}\big) = \notag\\
  &\min_{a,b} \ \ \sum_{j \in J^G} C_j g_j + \rho_{\xi^{t+1}}\left[{Q}_{t+1}\big(\{{v}^{t}_j, {a}^{[t-1]}_j\}_{j\in J^H}, \xi_{t+1}\big) \right]\label{mod:ht:init}\\
  &s.t.   \notag  \\
  & \sum_{j\in J^G_n} g_j + \sum_{j\in J^H_n} p_j u_j + \sum_{j\in J^R_n} r_j + \sum_{(n,m) \in J^F_n} f_{n,m} - f_{m,n}
  = {\tilde{D}_n}(\xi_t)  \label{mod:ht:load} \\
  & v^{t+1}_j = v^t_j - u_j - z_j + \sum_{n\in J^U_j} (u_n + z_n) + a_j^t, \quad \forall j \in J^H \label{mod:ht:mass} \\
  & 0 \leq v_j \leq V_j, \ \ 0 \leq u_j \leq U_j, \ \ 0 \leq z_j, \quad \forall j \in J^H \label{mod:ht:hyd} \\
  & 0 \leq g_j \leq G_j, \quad j \in J^G\label{mod:ht:ter}\\
  & 0 \leq r_j \leq \tilde{R}_j(\xi_{t}), \quad \forall j \in J^R \label{mod:ht:ren} \\
  & 0 \leq f_{n,m} \leq F_{n,m}, \quad \forall (n,m) \in J^F \label{mod:ht:fl} \\
  & a_j^t = \sum_{l=1}^L \phi_{j,l} a_j^{t-l} + \tilde{\varepsilon}_j(\xi_{t}), \quad \forall j \in J^H\label{mod:ht:end}
\end{align}

In the above model,
\eqref{mod:ht:init} is the objective function defined as the sum of the immediate cost (thermal generation cost)
plus the future cost.
\eqref{mod:ht:load} is the energy balance equation stating that the sum of energy generated in a node, plus the incoming energy minus the outgoing energy, must match the demand,
we assume there is enough thermal energy to match the demand in all stages and scenarios.
\eqref{mod:ht:mass} is the water mass balance stating that the reservoir level at the end of stage $t$ must match the initial content
minus the outflow (spilled and through turbine) plus the outflow of upstream plants plus the lateral inflow.
\eqref{mod:ht:hyd}-\eqref{mod:ht:fl} define physical bounds on variables.
\eqref{mod:ht:end} is the autoregressive equation of the inflow process that allows representing some non-stagewise independent processes in SDDP.

\subsection{Full tree simple example}

We first consider, as a simple example, a problem small enough
that we can compute its true optimal value by traversing
the full decision tree.
Our purpose is to view it as a sanity check to see that our upper bounds
are not only above the lower bounds but indeed above the true value.

We begin by considering a case with 10 stages,
2 scenarios per stage and, consequently, 512 possible scenario paths in the tree,
since the decision tree has exactly $\text{scenario}^{\text{stages}-1}$
possible paths.
We also compare the methods in the same case but
but altering the stages, scenarios and paths in the tree
to, respectively, 7, 3 and 729.
We show the results in table~\ref{tab:det-eqv}.

As one can see, the traditional sampling estimates a naive upper bound
that is consistently below the real cost calculated using the deterministic equivalent.
Meanwhile, the multicut-based sampling estimates approximately
the same cost as the full tree deterministic equivalent.
In the 2 scenarios cases, it is slightly above and in the 3 scenarios case
it is the same value.
This may happen because of round-off errors.

\begin{table}[thb]
    \caption{Total cost calculated via different methods.}
    \label{tab:det-eqv}
    \begin{center}
        \begin{tabular}{|c|c|c|}
            \hline
            Method & 2 scenarios & 3 scenarios \\
            \hline
            Deterministic Equivalent         & 685.4 & 462.4 \\
            Traditional Naive Upper Bound    & 636.4 & 444.7 \\
            Multicut-based Upper Bound       & 685.8 & 462.4 \\
            \hline
        \end{tabular}
     \end{center}
 \end{table}

\subsection{Brazil}

\begin{figure}[tbh]
  \centering
  \includegraphics[width=0.98\textwidth]{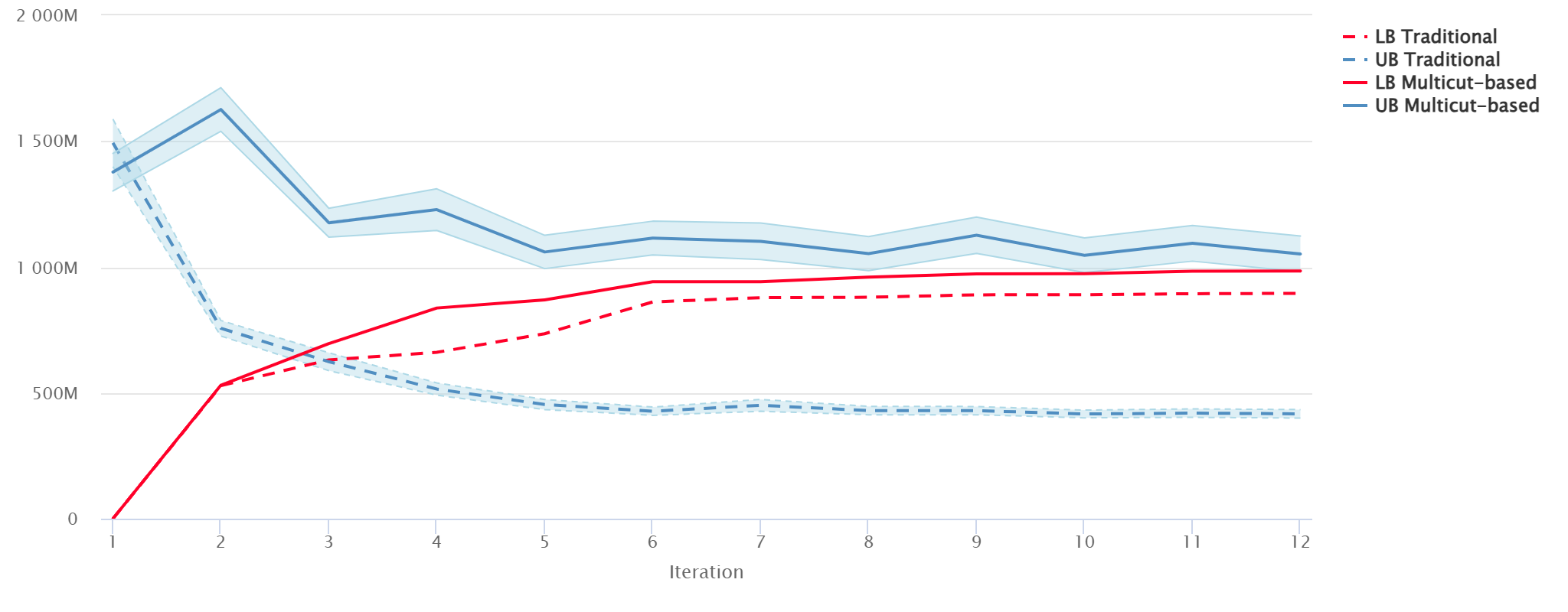}
  \caption{Convergence chart for the Brazilian system.}
  \label{fig:brazil}
\end{figure}

We compare the methods in the case study of~\cite{shapiro_case_br}
for the Brazilian interconnected power system operation planning.
The convergence charts for each method are shown in figure~\ref{fig:brazil}.

We see how the multicut-based sampling estimates for each iteration an upper bound above the lower bound
while the naive upper bound (known to be incorrect) computed with the traditional sampling does not.
It is also very interesting to notice that there is a meaningful increase in the calculated lower bound.
We attribute the improvement in the lower bound to an exploration of the scenario tree that is better aligned
with the nested CVaR objective.

\FloatBarrier

\subsection{Colombia}

Finally, we compare both sampling schemes for a case study
with real data from Colombia, including renewable plants and an autoregressive model for inflows.
It uses an interconnected network, and there is a non-trivial hydro topology,
as opposed to the other case studies in which hydro plants are in parallel.
The parameters are summarized in table~\ref{tab:parameters}.

Here we see the same results as before,
the multicut-based sampling estimates much better upper bounds
and also shows an increase in the value of the lower bound,
in comparison to the traditional sampling method.
We note here that the upper bound in the first iteration is below the lower bound of the second iteration.
We attribute this to the fact that, in the first iteration of SDDP, little is known about the future cost functions,
hence, the weights extracted from the optimal solution of the CVaR linear programming formulation might be
far from the optimal ones.

\begin{figure}[tbh]
  \centering
  \includegraphics[width=0.98\textwidth]{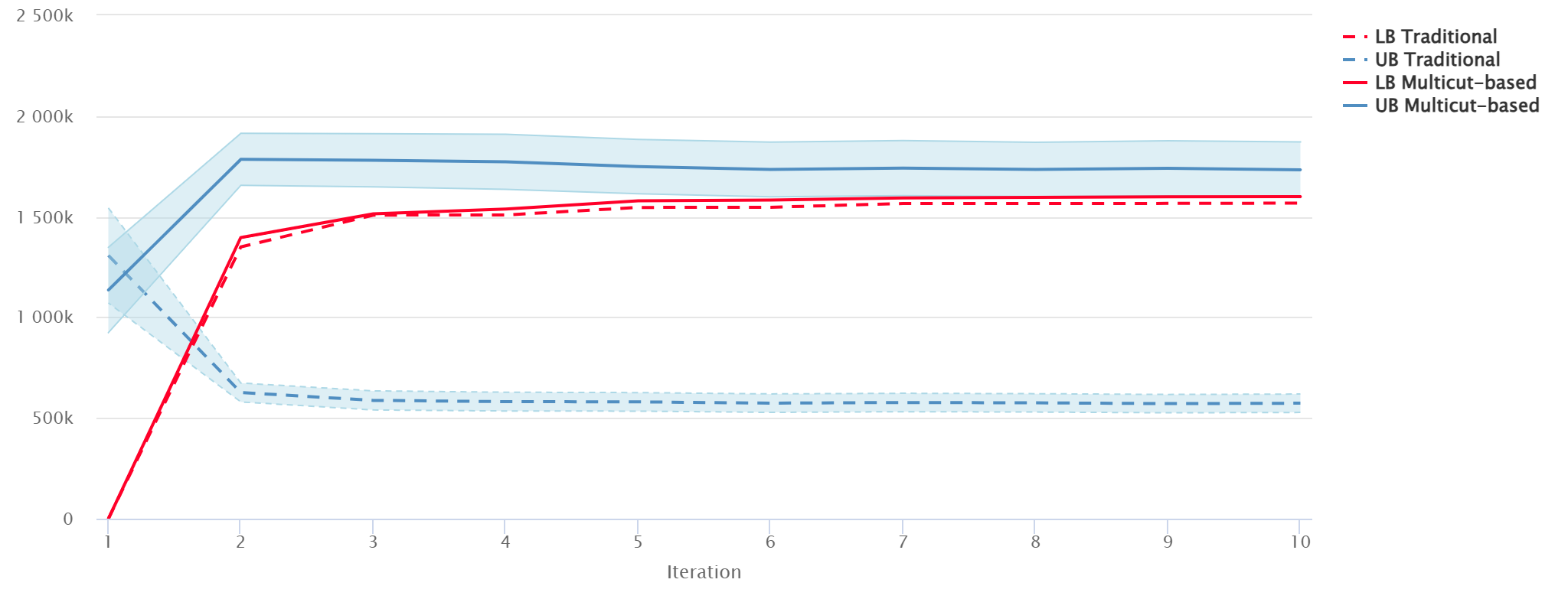}
  \caption{Convergence chart for the Colombian system.}
  \label{fig:colombia}
\end{figure}

\FloatBarrier

\section{Conclusions}

In this work, we presented a simple strategy to estimate the upper bound of a risk-averse multistage stochastic program solved with SDDP.
The methodology relies on the multicut reformulation to obtain weights used to sample scenarios of the next stage in the forward step of SDDP.
Such weights appear endogenously in the optimal solution of the CVaR linear formulation.
This methodology can be applied in the traditional convergence check of SDDP first proposed by \cite{sddp} and by other sample-based convergence criteria.

The simplicity of the methodology allows it to be easily implemented in existing SDDP code without major work.
Moreover, since the only change in the algorithm is the sampling in the forward pass, there is no performance degradation in each iteration.
Therefore, the methodology is computationally efficient.

We empirically demonstrated the effectiveness of this methodology by solving the full tree deterministic equivalent of a small example case and obtaining the same solution with the modified SDDP.
Moreover, we solved large hydrothermal dispatch problems: a standard academic version of the Brazilian system and a realistic representation of the Colombian system.
The solution of the latter two cases showed that the upper bounds behave much better than the naive approach and that the lower bounds were consistently improved.
We attribute the improved lower bounds to a better exploration of the scenario tree.

One possible drawback is that the methodology might be optimistic. The optimal resampling weights are only known when the Bellman functions are well approximated.
Otherwise, the sampling weights will be optimistic.

Finally, we highlight that there is a possible variant of this scheme for the single-cut case. The single-cut version of SDDP will not include the CVaR explicit linear formulation and hence, will not be able to obtain the weights directly from the solutions of the current linear program.
However, there is some information that could be cached while computing the single cuts, which is the weight used to average the cuts coming from each scenario.
If we cache this information, then after solving a subproblem combine the average of these weights according to the active cuts in the optimal solution, we can obtain weights to be used in the forward pass.
The main drawback of this version is that old cuts might have very poor weights associated with them, as they were obtained early in the algorithm.
Moreover, there is a need to cache more information associated to cuts.



\bibliography{refs}{}
\bibliographystyle{plain}

\end{document}